\documentclass[12pt]{article}

      \usepackage{latexsym}
         \usepackage[reqno, namelimits, sumlimits]{amsmath} 
         \usepackage{amssymb, amsfonts}
         \usepackage{amsthm} 

\newtheorem{Theorem}{Theorem}[section]

\newtheorem{theorem}[Theorem]{Theorem}
\newtheorem{Remark}[Theorem]{Remark}
\newtheorem{Remarks}[Theorem]{Remarks}
\newtheorem{lemma}[Theorem]{Lemma}
\newtheorem{proposition}[Theorem]{Proposition}

\newcommand{\vf}{\ensuremath{\varphi}}
\newcommand{\woo}{\ensuremath{W^{1,1}_0}}

\newcommand{\zr}{\ensuremath{Z_{\rm reg}}}
\newcommand{\zs}{\ensuremath{Z_{\rm sing}}}
\newcommand{\erf}{\ensuremath{{\rm erf}\,}}
\newcommand{\rt}{\ensuremath{\R\times(0,\infty)}}
\newcommand{\xoto}{\ensuremath{(x_0,t_0)}}
\newcommand{\xst}{\ensuremath{{{x}\over{\sqrt{2t}}}}}
\newcommand{\ist}{\ensuremath{{{1}\over{\sqrt{2t}}}}}
\newcommand{\sign}{\ensuremath{{\rm sign}\,}}
\newcommand{\lj}{\ensuremath{L^1(\R)}}
\newcommand{\ya}{\ensuremath{y_\alpha}}
\newcommand{\yast}{\ensuremath{\ya\sqrt{2t}}}
\newcommand{\heatdom}{\ensuremath{\R\times (0,\infty)}}
\newcommand{\ir}{\ensuremath{\int_\R}}
\newcommand{\tx}{\ensuremath{\tilde x}}
\newcommand{\ttt}{\ensuremath{\tilde t}}
\newcommand{\tw}{\ensuremath{\tilde w}}

\newenvironment{myproof}{\begin{proof}[{\rm\bf Proof.}]}{\end{proof}}

\renewcommand{\Re}{\operatorname{Re}}
\renewcommand{\Im}{\operatorname{Im}}

\newcommand{\R}{{\mathbb R}}

\newcommand{\be}{\beta}

\newcommand{\ep}{\epsilon}

\title{
Zeros of complex caloric functions and singularities
of
complex viscous Burgers equation}

 \author{P. Pol\'a\v cik\footnote{Supported in part by NSF Grant DMS-0400702}   
\ and 
V. \v Sver\'ak\footnote{Supported in part by NSF
  Grant DMS-0457061} 
\\
 School of Mathematics, University of Minnesota\\
 Minneapolis, MN 55455
}

\begin{document}

\numberwithin{equation}{section}

\maketitle

\begin{abstract}We show that the viscous Burgers equation
$u_t+uu_x=u_{xx}$ considered for complex valued functions
$u$ develops finite-time singularities from compactly supported 
smooth data. By means of the Cole-Hopf transformation,
the singularities of $u$ are related to zeros of
complex-valued  solutions $v$ of the heat equation $v_t=v_{xx}$.
We prove that such zeros are isolated if they are not present in the
initial data.
\end{abstract}

\section{Introduction}\label{sect1}
In a recent paper Sinai and Li \cite{Sinai} consider the
initial-value problem for the three-dimensional incompressible
Navier-Stokes equation, allowing the velocity field and the
pressure to be complex-valued. They prove that, in this setting,
there exist well-behaved (complex-valued) initial data for which
the solution blows up in finite time. In this note we consider a
similar problem for the viscous 1D Burgers equation
\begin{equation}
\label{B}
u_t+uu_x=u_{xx}
\end{equation}
in $\R\times(0,\infty)$ with initial condition $u(x,0)=u_0(x)$,
where we allow $u_0$ to be complex-valued. A well-known fact about
equation \eqref{B} is that the transformation $u=-2v_x/v$, called the
Cole-Hopf transformation, leads to standard heat equation
$v_t=v_{xx}$ for $v$. 
The singularities of $u$ correspond to the zeros of $v$. For real
valued functions, $v$ cannot have zeros if they are not present in
$v(x,0)$ and one sees immediately that for $u_0$ real and ``sufficiently
regular'' the initial value problem for equation~\eqref{B} has a
unique smooth global solution (in some natural classes of functions), see
\cite{Hopf}.
This can, of course, also be seen without the use of the Cole-Hopf
transformation, in a number of ways, since the equation \eqref{B}
has a maximum principle and an energy estimate with respect to
which it is subcritical. (The non-trivial scaling invariance of
equation \eqref{B} is the same as for Navier-Stokes:
$u(x,t)\to\lambda u(\lambda x,\lambda^2 t)$.)

The maximum principle and the energy estimates are lost when we
pass to complex-valued functions. At the same time, existence
proofs based on perturbation theory and Picard iteration, such as 
in \cite{Kato} or \cite{Koch-Tataru},   work also
in the complex case, and it is therefore natural to expect that
the proofs of local well-posedness for Navier-Stokes in critical
(i.\ e.\ scale-invariant) spaces work also for equation~\eqref{B},
without using its ``complete integrability''. One can therefore
expect local well-posedness of complex-valued equation~\eqref{B} in
$L^1$ (by analogy with \cite{Kato}) and, in fact, in $(BMO)^{-1}$
(by analogy with \cite{Koch-Tataru}). With the Cole-Hopf
transformation, the local $L^1$ well-posedness becomes completely
transparent (see below). It it not quite so with the local
$(BMO)^{-1}$ well-posedness, which shows the subtle nature of the
well-posedness result in \cite{Koch-Tataru}. Our focus will be on
global well-posedness, and therefore we will work with the $L^1$
space which is very simple and - as we will see - completely
adequate for the problems we will consider.

Since the zeros of $v$ produce singularities of $u$, it is easy to
find compactly supported smooth (complex-valued)
initial data $u_0$ for which the solution of equation~\eqref{B}
blows up in finite time, see Proposition~\ref{blowup}. 
The continuity argument used in the proof of this proposition
allows us to formulate more general sufficient conditions on $u_0$ 
for the solution $u$ to blow up,  see Remark~\ref{rm-sing}(i). 
On the other hand, in Proposition~\ref{glex} we give a sufficient 
condition on $u_0$ for the solution to converge to zero. Using these
results, we can then explicitly  describe the boundary,
in some subsets of the space of initial data $L^1$,
between the basin of attraction of zero and the region from which the
solutions blow-up, see Remark~\ref{rm-exp}(iii). 
The behavior of solutions with the initial conditions 
on this boundary is then naturally of interest. These solutions 
are global and bounded and we describe their asymptotics 
as $t\to\infty$, see Proposition~\ref{asympt} and Remark~\ref{rm-exp}(iii).

Once we know that a solution can develop singularities,
we can  ask about the nature of the
singular set. We will prove that, roughly speaking, if there are
no singularities present in the initial data, then the set of
singularities of the function $u$ defined by the Cole-Hopf
transformation from $v$ (that is, $u=-2v_x/v$) is always discrete
in $\R\times (0,\infty)$.
This follows from a theorem about zeros of complex-valued
solutions of 1d heat equation (Theorem~\ref{zeros}). 
In a ``typical situation'' the number of singularities of such a
solution $u$ will be finite. However, as we  show in Section~\ref{infinite},  
certain regular initial data yield solutions 
with infinitely many (isolated) singularities. 
We will also briefly address the question what ``right-hand side'' 
the singularities produce in a suitable weak formulation of the
equation (see Section \ref{additional}).

The solution of equation~\eqref{B} defined by $u=-2v_x/v$ is
analytic outside a discrete set. This is not a typical behavior of
solutions of non-linear parabolic equations with singularities. In
fact, it is reasonable to expect that, for many equations, analyticity
in the time variable will be destroyed in the whole time level
$\{(x,t), t=t_0\}$ if we have a singularity at time $t_0$ at some
point $x_0$. This conjecture is based on the study of singularities
of the Complex Ginzburg-Landau equation in  \cite{plsv}.
As far as we know, the issue has not been much
studied.

In the case of the dispersive regularization of Burgers equation,
which is the KdV equation $u_t+uu_x=u_{xxx}$, the singularities
for complex-valued solutions are studied in~\cite{birnir}.
Viscous Burgers equation with complex viscosity is studied
by means of the Cole-Hopf transformation in~\cite{senouf1}
and~\cite{senouf2} for a particular real initial condition,
with the main focus on the behavior of singularities
arising in complex time.

\section{Cole-Hopf transformation and singularities}\label{sect2}
For a complex-valued $u\in L^1(\R)$ we define
$U(x)=\int_{-\infty}^x u(\xi)\,d\xi$ and $v=\exp(-U/2)$.
Vice-versa, given a complex-valued $v\in\woo(\R)$ (the space of
all absolutely continuous functions that have  the derivative in 
$L^1(\R)$) with $v(x)\ne0$
in $\R$ and $v(x)\to 1$ as $x\to -\infty$, we let $u=-2v_x/v$. For
time-dependent functions on $\R$ we apply the above transformations
at each time level.

A well known simple calculation shows that $u$ satisfies
equation~\eqref{B} if and only if $v$ satisfies the standard heat
equation $v_t=v_{xx}$, see for example \cite{Hopf}. (If one does
not impose the normalization $v(x)\to 1$ as $x\to -\infty$, the
function $v$ is only determined up to a multiplicative factor  
depending on time, and the heat equation for $v$ needs an extra term
which would account for this, see \cite{Hopf}.)

We can now solve the initial value problem for equation~\eqref{B}
with a complex valued $u_0\in L^1(\R)$ as follows. Set
$v_0(x)=\exp\{-{1\over2}\int_{-\infty}^x u_0(\xi)\, d\xi\}$ and
let $v$ be the bounded solution of the heat equation with initial
data $v_0$. It is easy to check that there is $T>0$ such that
$|v(x,t)|\ge\varepsilon>0$ in $\R\times(0,T)$ and hence $u =
-2v_{x}/v$ is a well-defined local-in-time solution of
equation~\eqref{B} with $u(x,0)=u_0(x)$.

\begin{proposition}
\label{glex} In the notation above, assume that $u_0\in L^1(\R)$ with
$\int_\R|\Im u_0|\le 2\pi$. Then equation~\eqref{B} has a global smooth
solution $u$ with $u(x,0)=u_0(x)$. If in addition
$|\int_\R\Im u_0|<2\pi$, then
 $\,\sup_x|u(x,t)|\to 0$ as
$t\to\infty$.
\end{proposition}

\begin{myproof}
 When $\int_\R|\Im u_0|\le 2\pi$, the function
$v_0(x)=\exp\{-{1\over2}\int_{-\infty}^x u_0(\xi)\, d\xi\}$ takes
 values in a convex sector of the form $\{z,\,\alpha<\arg(z)<\beta\}$ with
$\be-\alpha\le\pi$ and it has finite nonzero limits 
as $x\to\pm\infty$. Thus for suitable real $\theta$ we have
 $\Re\, (e^{i\theta} v_0(x)) >0, \, x\in \R$. Applying the 
strong maximum principle  to $\Re\,  (e^{i\theta}v)$ (which  solves 
the heat equation), we see that  $\Re\,  (e^{i\theta}v)>0$. Therefore
 $v$ cannot vanish
at any point $(x,t)\in \R\times(0,\infty)$, 
proving the first statement. 
If $|\int_\R\Im u_0|<2\pi$, we can choose $\theta$ so that
$\Re (e^{i\theta} v_0(\pm\infty))\ge \varepsilon_1>0$, and therefore 
for all large
$t$ we will have $|v(x,t)|>\varepsilon_1/2$.
Moreover, since $v_x$ also solves 
the heat equation and $v_{0,x}\in L^1(\R)$, we have 
$\,\sup_x|v_x(x,t)|\to 0$ as
$t\to\infty$. These  properties imply
the second statement.
\end{myproof}

\begin{proposition}
\label{blowup}
For each $\delta>0$ there exists a smooth, compactly supported
(complex-valued) $u_0$ with $\int_\R |u_0|<2\pi+\delta$ such that the
solution of equation~\eqref{B} with initial condition $u_0$ blows up
in finite time.
\end{proposition}

\begin{myproof} We choose a smooth compactly supported non-negative $\vf$
with $\int_\R\vf=2\pi+\delta/2$. Set $u_0=-i\vf$ and let $v_0$ be
the Cole-Hopf transformation of $u_0$. The function $v_0$
satisfies:
\begin{itemize}
\item $v_0(x)=1$ for large negative $x$,
 \item
$v_0(x)=\exp(i(\pi+\delta/4))$ for large positive $x$,
 \item
$0\le\arg(v_0(x))\le\pi+\delta/4$ for $x\in \R$.
\end{itemize}
Let $v$ be the solution of the heat
equation with initial data $v_0$. For $\theta\in \R$ 
the function $e^{-i\theta} v$ also solves 
the heat equation and choosing $\theta>0$
small enough we achieve that 
\begin{equation*}
  \lim_{x\to\pm\infty}\Im (e^{-i\theta}v_0(x))<0. 
\end{equation*}
It then follows that for a sufficiently large  $t_0>0$ we have
$\Im (e^{-i\theta}v(x,t_0))<0$ for all $x\in\R$. Since
the limits of $e^{-i\theta}v(x,t)$ as $x\to\pm\infty$ are
independent of $t$, 
comparing the trajectories of $x\mapsto e^{-i\theta}v(x,t)$ for 
$t=0$ and $t=t_0$  we conclude that $v$  
has to vanish at some point $(x_1,t_1)$ with $t_1\in (0,t_0)$.
Consequently $u$ has a singularity at $(x_1,t_1)$.
\end{myproof}

\begin{Remarks} {\rm
  \label{rm-sing}
(i)   
 The above continuity argument can also be used  to show that 
$u$ develops a singularity whenever $u_0\in L^1(\R)$ 
satisfies 
 $ |\int_\R \Im u_0|> 2\pi$ and $\int_\R \Im u_0$ is {\it not}
 of the form $2\pi+4k\pi$, where $k$ is an integer.
 
(ii) It is perhaps worth pointing out a  very non-local behavior
implied by Propositions~\ref{glex}, \ref{blowup}. 
Consider a compactly supported
$u_0$ with $\Im u_0\ge 0$ and $\ir \Im u_0=\pi+\delta$ for some 
small $\delta >0$. Then the solution of equation~(\ref{B}) with initial
condition $u_0$ exists for all time and converges to zero.
Consider now the initial condition
$\tilde u_0^a(x)=u_0(x-a)+u_0(x+a)$. With initial condition $\tilde u^a_0$,
the solution of (\ref{B}) will blow up, no matter how large $a$ is.
If we take $a$ very large, the solution will become very small in
$L^\infty$ before it starts growing again and blows up.
(In fact, it is not hard to see that one can replace 
$L^\infty$ by $L^p$ for a fixed
$p>1$ in the last sentence.)
}\end{Remarks}

\begin{proposition}
\label{asympt}
Assume $u_0\in \lj$ is compactly supported, with 
$|\int_\R \Im u_0|=\int_\R |\Im u_0|=2\pi.$ Then there exist
a real $\ya$ and a complex $\beta$ with $\Im \beta\ne 0$
such that the solution $u$ of equation~\eqref{B} with
$u(x,0)=u_0(x)$ satisfies 
\begin{equation}
\label{asymptotics}
u(x,t)={{-2}\over{(x-\yast)+\beta}}\,\,+\,O({1\over{\sqrt{t}}}\,), \quad (t\to\infty)\,
\end{equation}
uniformly in $x\in\R$.
\end{proposition}

\begin{myproof} By Proposition~\ref{glex}, the solution $u$ is global.
Let $U_0(x)=\int_{-\infty}^xu_0(\xi)\,d\xi$ and $ v_0(x)=\exp(-U_0(x)/2)$.
Let $L>0$ be such that $u_0$ vanishes outside $[-L,L]$. We have $v_0(x)=1$ on
$(-\infty, -L)$ and $v_0(x)=-\exp(-I/2)$ on $(L,\infty)$, where $I=\int_R \Re u_0$.
We set
$F(x)=(2\pi)^{-1/2}\int_0^x \exp(-\xi^2/2)\,d\xi$.
(In terms of the $\erf$ function used in probability we can
 write $F(x)=\erf(x)-1/2$.) We note that the function $F(\xst)$
 solves the heat equation in $\R\times (0,\infty)$ with 
 the initial data ${1\over2}\,\sign(x)$. The fundamental solution
 $\Gamma(x,t)$ of the heat equation can be written as 
 $\Gamma(x,t)=\ist F'(\xst)$.
 We write the initial condition $v_0$ in the form
 \begin{equation}
 \label{initial}
 v_0(x)=-a\,\frac{\sign(x)}{2}+b+w_0(x)
 \end{equation}
 where $a,b$ are chosen so that $w_0$ be compactly supported.
 This gives $a=1+\exp(-I/2)$ and $b=(1-\exp(-I/2))/2$. Also note
 that  $\Im w_0$ is continuous and does not change sign.
 Let $w(x,t)$ be the solution of the heat equation in $\heatdom$
 with initial condition $w_0$. The solution with the initial condition
 $v_0$ is then
 \begin{equation}
 \label{vform}
 v(x,t)=-aF(\xst)+b+w(x,t)\,.
 \end{equation}
 From the  representation formula $w(x,t)=\ir w_0(y)\Gamma(x-y,t)\,dy$
 we see that 
 \begin{equation}
 \label{decayw}
 |w(x,t)|+\sqrt{t}\,|w_x(x,t)|=O({1\over\sqrt{t}}), \, \quad (t\to\infty)
 \end{equation}
 uniformly in $x$.
 To get further estimates for $w$ in a simple way, 
 we will use Appell's transformation and write
 \begin{equation}
 \label{appell}
 w(x,t)=\Gamma(x,t)\,\tilde w({x\over t},-{1\over t})\,,
 \end{equation}
 where $\tilde w$ is a function of $\tilde x=x/t$ and $\tilde t=-1/t$ (defined 
 for $(\tilde x,\tilde t)\in \R\times(-\infty,0)$) which again satisfies
 the heat equation in $\tx$ and $\ttt$. 
 We have
 \begin{equation}
 \label{bounded}
 \tw(\tx,\ttt)=\ir w_0(y){{\Gamma(x-y,t)}\over{\Gamma(x,t)}}\,dy=
 \ir w_0(y)\exp(-{{\tx y}\over2}+{{\ttt y^2}\over4})\,dy\,,
 \end{equation}
 which shows that $\tw$ can be analytically extended to $\R\times\R$.
  Letting $c=\ir w_0$, we see from~(\ref{bounded}) that $\tw(0,0)=c$.
 Setting $z=\tw-c$, we can write
 \begin{equation}
 \label{repres}
 v(x,t)=-a\,F(\xst)+b+c\,\Gamma(x,t)+z({x\over t},-{1\over t})\,\Gamma(x,t)\,,
 \end{equation}
 where $z$ is smooth in $\R\times(-\infty,0]$ and $z(0,0)=0$.
 We now set $\alpha=b/a$ and $\beta=-c/a$. Observe that
 $\alpha\in(-1/2,1/2)$ and $\Im\beta\ne 0$. 
From~(\ref{repres}) we get
 \begin{equation}
 \label{equ}
 u=-2{\frac{v_x}{v}}=-2 \frac{
{1+\beta{ \frac{\Gamma_x}
 {\Gamma}}
 -{\frac{z}{a}}{\frac{\Gamma_x}\Gamma} 
 -{\frac{z_{\tx}}a}{\frac{1} t}}}
 {\frac{F-\alpha}{\Gamma}+\beta-\frac{z}{a}}\,,
 \end{equation}
 where $F$ is evaluated at $x/\sqrt{2t}$.
 Given any fixed $R>0$ and $\ttt_0<0$ we can see that 
 for $|\tx|<R, \,\,\ttt>\ttt_0$ we have
 $|z_{\tx}|\le c_1,\, |z|\le c_2|\tx|+c_3|\ttt|$. Together
 with~(\ref{equ}),  and after taking into account that $\Im \beta\ne
 0$, 
this gives
 \begin{equation}
 \label{equ2}
 u=-2{{1+O({1\over\sqrt{t}})}
 \over
 { {{F-\alpha}\over{\Gamma}}+\beta+O({1\over\sqrt{t}})}}
 ={{-2}\over{{F-\alpha}\over\Gamma}+\beta}+O({1\over\sqrt{t}}\,)\,,\quad 
(t\to\infty) 
 \end{equation}
 uniformly in regions $\{(x,t);\,|x|/\sqrt t\le R\}$, 
 where $F$ is again evaluated at $x/\sqrt{2t}$.
 
 Let $\ya$ be the unique root of the equation $F(y)=\alpha$ and
 let us fix a (small) $\delta_1>0$. We note that when
 $|x/\sqrt{2t}-\ya|\ge\delta_1$, then $|aF(x/\sqrt{2t})-b|\ge \varepsilon_1>0$,
 and from \eqref{vform}, \eqref{decayw} we see that, in the region
 $\{(x,t);\,|x/\sqrt{2t}\,\,-\ya|\ge\delta_1\}$, one has
 \begin{equation*}
u=-2{{v_x}\over{v}}=O({1\over{\sqrt{t}}}),\quad (t\to\infty), 
\end{equation*}
uniformly in $x$ (for $(x,t)$ in the above region).
Taking into account~(\ref{equ2}), we  see that it only remains
to show that 
\begin{equation}
\label{remains}
{{-2}\over{{F-\alpha}\over\Gamma}+\beta}=
{{-2}\over{(x-\yast)+\beta}}\,\,+\,O({1\over{\sqrt{t}}}\,), \quad (t\to\infty)
\end{equation}
uniformly in $\{(x,t);\,|x/\sqrt{2t}-\ya|\le\delta_1\}$. (As above, $F$ is evaluated
at $x/\sqrt{2t}$.)
For $y\in \R$ we set
\begin{equation*}
\kappa(y)={{F(y)-F(\ya)}\over{F'(y)(y-\ya)}}\,
\end{equation*}
with the understanding that $\kappa(\ya)=1$.
Since $\Gamma(x,t)={1\over\sqrt{2t}}F'(\xst)$, we have
\begin{equation}\label{cor}
{{-2}\over{{F(\xst)-\alpha}\over\Gamma(x,t)}+\beta}=
{{-2}\over{\sqrt{2t}{{F(\xst)-F(\ya)}\over{F'(\xst)}}+\beta}}=
{{-2}\over{\kappa(\xst)(x-\ya\sqrt{2t})+\beta}}
\end{equation}
Clearly
\begin{equation}
\label{kappaineq}
\kappa(y)=1+O(|y-\ya|)\,,\quad (y\to\ya).
\end{equation}
Using the elementary formula
\begin{equation*}
{1\over{\kappa\xi+\beta}}-{1\over{\xi+\beta}}=
\int_0^1{{\xi(1-\kappa)}\over{(\xi[(1-s)+s\kappa]+\beta)^2}}\,ds
\end{equation*}
we see that~(\ref{remains}) follows from~(\ref{cor}) and (\ref{kappaineq}).
\end{myproof}

\begin{Remarks}
  \label{rm-exp}
 {\rm
(i) We note that, with the notation used in the proof, we have
    $\Im\beta={1\over a} \ir\Im w_0$.   Under the assumptions
    of Proposition~\ref{asympt} the function $\Im w_0$ does not
    change sign and is integrable. If we do not assume that
    $u_0$ is compactly supported, it can easily happen
    that $|\ir \Im w_0|=+\infty$. We then have $|\Im \beta|=+\infty$
     and in view of~(\ref{asymptotics}) 
    it is natural to expect that in that case $u(x,t)\to 0$ as 
    $t\to\infty$ uniformly in $x$. 
    This can indeed be proved. If $u_0$ is not compactly
    supported, but $\int \Im w_0$ is finite, we expect that the 
    asymptotics of $u(x,t)$ will be similar to~(\ref{asymptotics}),
    with perhaps a slower rate of convergence.

(ii) The constant $\ya$ is given by the equation 
$F(\ya)=\frac{1}{2} \tanh({I\over4})$, with
$I=\ir \Re u_0$. In particular, $\ya=0$ if and only if $\ir \Re u_0=0$. 
We note that ${{-2}\over{x+\beta}}$ 
is a steady-state solution of equation~(\ref{B}).
 For $I\ne0 $ we have 
 ${{d}\over{dt}}(\ya\sqrt{2t})={{y_\alpha}\over{\sqrt{2t}}}\ne 0$,
 and we can interpret formula~(\ref{asymptotics}) as an
 ``almost steady state solution'', which is slowly drifting
 to $\pm\infty$, at speed~${{y_\alpha}\over{\sqrt{2t}}}$.

 (iii) Consider a complex valued $L^1$ function $\tilde u_0$
 supported in $[-L,L]$, with $\int_\R |\Im \tilde u_0|=\int_\R\Im
 \tilde u_0=2\pi$.  
 Let $\mathcal O(\tilde u_0,L,\varepsilon)=\{u_0\in L^1(\R):\,
 ||u_0-\tilde u_0||_{L^1}<\varepsilon,\, 
 \mbox{$u_0$ is supported in $[-L,L]$}\}$.
 From the above one can see that for sufficiently small $\varepsilon$,
  one has, in the set 
  $\mathcal O(\tilde u_0,L,\varepsilon)$, an explicit description
 of the boundary between the basin of attraction of the zero
 solution of equation~(\ref{B}) and the region from which
 the solutions of equation~(\ref{B}) blow up is finite time:
The boundary (in $\mathcal O(\tilde u_0,L,\varepsilon)$) is given
by the equation $\int_\R \Im u_0 = 2\pi$.
(To be precise, for the proof of this one needs to augment the above
propositions by a slightly modified 
version of  Proposition~\ref{glex}, in which we assume that $u_0$ is in
$\mathcal O(\tilde u_0,L,\varepsilon)$, $\int_\R \Im u_0<2\pi$, and 
we allow
$\int_\R |\Im u_0|\le 2\pi+\delta$, where 
$\delta=\delta(L,\varepsilon,\tilde u_0)>0$ is sufficiently
small. The restriction on the support of $u_0$
is crucial in this step. We leave the details to the reader.)
It is not hard to check that the large-time asymptotics of the 
solutions starting at the boundary 
(in $\mathcal O(\tilde u_0,L,\varepsilon)$) of the basin of attraction
of the zero solution is given
by the solutions described in Proposition~\ref{asympt}. 
If we replace $\mathcal O(\tilde u_0,L,\varepsilon)$
by $\mathcal O(\tilde u_0,\infty,\varepsilon)$ (i.\ e.\ we remove
the restriction on the support of the perturbed function), the situation
changes and the boundary is no longer described in a simple way.
In addition, even when $\varepsilon$ is small,  we expect that for
$\mathcal O(\tilde u_0,\infty,\varepsilon)$ some solutions at the boundary 
of the basin of attraction of zero will have more
complicated behavior, such as slow oscillations with large
amplitude. 
}
\end{Remarks}

\section{Nodal sets of caloric functions}\label{set3}
Let $u$ be a bounded real-valued nontrivial 
solution of the heat equation
$u_t=u_{xx}$ in \rt. We define

 $\begin{array}{lll}
Z & = & \{(x,t)\in\rt;\,u(x,t)=0\}\,, \\
\zr & = & \{(x,t)\in Z;\,u_t^2+u_x^2\ne 0\}\,, \quad\mbox{and}  \\
\zs & =& Z\setminus\zr\,.
\end{array} $

\noindent
 The analyticity of $u$ implies that $\zs$ is discrete
and that $Z$ is locally a regular real analytic curve in a
neighborhood of each point $\xoto$ in $\zr$. By a 
regular (real) analytic curve $C$ 
in an open set $U\subset \rt$ we mean
 a one-dimensional analytic (imbedded) submanifold of 
$U$ with $\bar C\setminus C\subset \partial U$. 

\begin{lemma}
\label{continuation}
 In the notation introduced above, the regular analytic
curves describing $Z$ in a neighborhood of $\xoto\in\zr$ can be
analytically continued through the points of $\zs$. In other
words, $Z$ is a (locally finite) union of regular analytic curves
in $\rt$.
\end{lemma}

\begin{myproof}
We first recall some facts about caloric polynomials. As
usual in the parabolic setting we say that a function $f(x,t)$ is
parabolically $m-$homogeneous if $f(\lambda x,\lambda^2
t)=\lambda^m f(x,t)$ for $\lambda >0$. The $m-$th caloric
polynomial is a parabolically $m-$homogeneous polynomial
satisfying the heat equation. It is unique, modulo a
multiplicative factor, and can be given for example by
\begin{equation}
\label{hermite}
P_m(x,t)=\sum_{k=0}^{[m/2]}{{m!}\over{k!(m-2k)!}}\,x^{m-2k}t^k.
\end{equation}
The polynomial $P_m(x,-1)$ is the $m-$th Hermite polynomial. For $m$ even, 
$m=2k$, the polynomial $P_m$ is of the form
\begin{equation}
\label{even}
P_m(x,t)=(x^2+a_1t)\cdots(x^2+a_kt),
\end{equation}
with $0<a_1<\cdots <a_k$. For $m$ odd, $m=2k+1$, $P_m$ is of the
form
\begin{equation}
\label{odd}
P_m(x,t)=x(x^2+a_1t)\cdots(x^2+a_kt),
\end{equation}
with $0<a_1<\cdots <a_k$. (The $a_j$'s may be different for
different $m$, of course.)

In a neighborhood of a point $\xoto\in\zs$ we can write $u$ as a
convergent series
\begin{equation}
\label{series} 
u(x,t)=a_m P_m(x-x_0,t-t_0)+a_{m+1}P_{m+1}(x-x_0,t-t_0)+\cdots
\end{equation}
where $m\ge3,\,a_m\ne0$. In what follows 
we will assume that $m$ is odd, $m=2k+1$. (The proof for $m$ even is similar 
and, in fact, easier.)
We change coordinates so that $\xoto$ corresponds to $(0,0)$ in the new 
coordinates, which we still denote $(x,t)$. We let $t=-y^2, a_1=b_1^2,
\cdots,a_k=b_k^2,\,\,b_j>0$. The equation $u(x,t)=0$ can be written as
\begin{equation}
\label{eqxy}
u(x,-y^2)=x(x-b_1y)(x+b_1y)\cdots(x-b_ky)(x+b_ky)+R(x,y)=0,
\end{equation}
where $R(x,y)$ is analytic in a neighborhood of $(0,0)$ with
vanishing derivatives of order $1,2,\dots m$. 
Letting $b_0=0$, we will look for analytic curves $x=x(y)$ of the
form $x(y)=by+y^2f(y)$ (with $b=\pm b_j,\,j=0,1,\dots,k$)  defined
for small $y$ on which $u(x,-y^2)$ vanishes. Substituting the expression
$x(y)=by+y^2f(y)$ in equation~\eqref{eqxy}, it is easy to check that
we get an equation of the form
\begin{equation}
\label{feq}
f(y)=F(y,f(y))
\end{equation}
where $F=F(y,f)$ is analytic in $f$ and depends on $f$ only
through $yf$. Therefore $f_0:=F(0,f)$ is independent of $f$
and ${{\partial F}\over{\partial f}}(0,f)=0$.
Applying the standard implicit function theorem one
shows that equation~\eqref{feq} has an analytic solution 
$f$ defined on a neighborhood of 0 with $f(0)=f_0$.

Observe that with $y\ne 0$ fixed, the curves found above give
$m$ different solutions of \eqref{eqxy}. Therefore, 
by the Weierstrass  preparation theorem, they yield all 
solutions of \eqref{eqxy} in a neighborhood of 
$(0,0)$. (One can also use the Malgrange preparation theorem.)

To finish the proof, we note that for $j\ge 1$, instead of writing
$x=x(y)$ we can write $y=y(x)$ and the equation $t=-y^2=-(y(x))^2$
then defines the analytic branch of $Z$ which has contact of the second 
order with the parabola $t=-x^2/a_j$. When $j=0$ we note that the 
function in equation~\eqref{feq} is of the 
form $F(y,f)=\tilde F(-y^2,f)$,
and hence the corresponding curve is of the form
$x=x(t)=t\tilde f(t)$, with an analytic $\tilde f$ satisfying  $\tilde
f(0)\ne 0$.
\end{myproof}

\begin{Remark}
  \label{rm-nodal}
{\rm
It is clear that the proof of the lemma also works
when the equation has lower-order terms and analytic
coefficients.
Although we did not find the precise
statement of the lemma in the literature,   
we assume it is  known to experts.  For example,
it follows easily from the analysis of nodal sets
in \cite{Angenent-F:s1}, where the method of Newton
polygons is used. Also, once the specific form of the 
caloric polynomials
is taken into account, the lemma can be derived 
easily from general principles used
in algebraic geometry for ``desingularization''.
 Nevertheless,
we think that the elementary proof above is still 
of some interest and we have included it for completeness.
We remark that even if we allow nonanalytic variable coefficients, 
$Z$ is  still a finite union
of regular $C^1$ curves in a neighborhood of any point in 
$Z_{\rm sing}$, see \cite{Chen:strong}.  
}
\end{Remark}

\begin{theorem}
\label{zeros}
Let $v$ be a bounded complex-valued solution of the heat equation in $\rt$.
Assume $\,v$ has no zeros in some neighborhood of  
$\R\times\{0\}$. Then all zeros
of $\,v$ in $\rt$ are isolated.
\end{theorem}

\begin{myproof}
Let $v=v_1+iv_2,\,\,Z_1=\{v_1=0\},\,\,Z_2=\{v_2=0\}$. Assume $Z_1\cap Z_2$
has an accumulation point $\xoto$ inside $\rt$. By Lemma~\ref{continuation}
we know that $Z_1$ is a locally finite union of regular analytic
curves in $\rt$. 
Consider the curves passing through $\xoto$. Clearly
$v_2$ has infinitely many zeros 
accumulating at $\xoto$ on one of the curves, let's call it $C$.
By analyticity, $v_2$ vanishes on $C$. The
curve $C$ cannot be closed, for otherwise the maximum principle
would imply that both $v_1$ and $v_2$ 
vanish in the interior of $C$, which is impossible
by our assumption and analyticity. 
Hence we can parametrize $C$ by a parameter $s\in(-\infty,\infty)$. Also, 
since $C$ is a regular analytic curve, we have 
$\bar C\setminus C\subset\R\times\{0\}$, thus our assumption
implies that $\bar C\setminus C=\emptyset$.
Now either  the time coordinate $t$ has a strict local minimum on $C$
or we can choose the parametrization so that $t(s)$ is monotone 
nonincreasing  for large $s$ and $x(s)$ approaches 
$\infty$ or $-\infty$
as $s\to\infty$. In either case, we find a (bounded or unbounded)
domain in $\rt$ such that both 
functions $v_1$ and $v_2$ vanish on its parabolic boundary.
Since they are bounded, the maximum principle \cite{Friedman:bk-par} 
implies that they vanish on a nonempty open set, 
hence on $\rt$, and we again have a  contradiction to our assumption.

\end{myproof}

 It is clear that the proof of Theorem~\ref{zeros} works without 
much change also for complex-valued harmonic functions in a half-plane.

\section{Additional comments on the singularities}\label{additional}
Given complex-valued initial data $u_0\in L^1(\R)$ for equation~\ref{B}
and constructing  
the solution $u$ of the initial-value problem by means of the
Cole-Hopf transformation 
as in Section~2 by setting $u=-2v_x/v$, 
we see from Theorem~\ref{zeros} that the
singularities of $u$ are isolated. It is natural to ask if
equation~\eqref{B} is satisfied 
in some weak sense across the singularities, or if the singularities
introduce a non-trivial 
``right-hand side'', i.\ e.\ we want to calculate the distribution $f$ given by
$u_t+uu_x-u_{xx}=f$, where the left-hand side requires a suitable interpretation.
Clearly $f$ should be supported in the singular set.
Even if we write the operator $u_t+uu_x-u_{xx}$ as 
$u_t+(u^2/2)_x-u_{xx}$ the definition of $f$  can still be somewhat
ambiguous, since 
$u$ and $u^2$ are not locally integrable in a neighborhood of a singularity.
We suggest one possible interpretation. For simplicity we will consider only
the simplest case when the function $v$ defining $u$ has a simple zero
at the singularity  
$\xoto$, i.\ e.\
$v(x,t)=a(x-x_0)+b(t-t_0) + O((x-x_0)^2+(t-t_0)^2)$, with $a,b$ complex
and linearly  
independent over $\R$.  We consider a smooth test function
$\vf=\vf(x,t)$ supported 
in a small neighborhood of $\xoto$, so that no other singularity is
present in the 
support of $\vf$. We want to define 
\begin{equation}
\label{integral}
I=\int_{\rt}(-u\vf_t-u^2{\vf}_x/2-u{\vf}_{xx})\,dx\,dt. 
\end{equation}
For $t\ne t_0$ we let 
\begin{equation}
\label{h}
h(t)=\int_\R(-u(x,t)\vf_t(x,t)-u^2(x,t){\vf}_x(x,t)/2-u(x,t){\vf}_{xx}(x,t))\,dx
\end{equation}
Using well-known facts about the behavior of the 
distributions $1/(x\pm i\varepsilon)$ 
and $1/(x\pm i\varepsilon)^2$ as $\varepsilon\to 0$, together with a
change of variables 
$v(x,t)=y$ (for a fixed $t$) one can see that $h(t)$ has one-sided
limits as $t\to t_0$. 
Therefore it seems to be natural to define integral $I$ in
expression~\eqref{integral}  as 
\begin{equation}
\label{hh}
I=\int_0^\infty h(t)\,dt=\lim_{\tau\to 0}\big(\int_0^{t_0-\tau}h(t) 
\,dt+\int_{t_0+\tau}^\infty h(t)\,dt.\big)
\end{equation} 
Integration by parts, the equation satisfied by $u$, 
and the specific form of the singularity of $u$ now give
\begin{equation}
I=\lim_{\tau\to 0} \int_\R\big[u(x,t_0+\tau)\vf(x,t_0+\tau)-
u(x,t_0-\tau)\vf(x,t_0-\tau)\big]\, dx = \pm4\pi i\vf(x_0),
\end{equation}
where the correct sign is the same as the sign of the imaginary part of
$\bar a b$. 
In other words, in a neighborhood of $\xoto$  we have, in some sense,
\begin{equation}
u_t+uu_x-u_{xx}=\pm4\pi i\delta_{\xoto},
\end{equation}
where we use the usual notation $\delta_{\xoto}$ for the Dirac distribution
at the point $\xoto$.
Therefore we cannot interpret $u$ as a global weak solution of
equation~\eqref{B}.

\section{Infinitely many singularities}\label{infinite}
If is not hard to show that, typically, the solution 
$u$  will only have finitely many singularities. 
In fact, one can check easily that a sufficient condition
for $u$ to  have only
finitely many singularities is that
 $\int_\R \Im u_0(x)\,dx$ is {\it not} of the form
$2\pi+4k\pi$ with $k$ an integer. 
We now show that, on the other hand, there are solutions with
regular initial data having 
infinitely many singularities. This is an immediate
 consequence of Proposition~\ref{propin} below.
We recall that we denote by $\woo(\R)$ the space of all
functions on $\R$ which are absolutely continuous with
the derivative in $L^1(\R)$. (In particular,  
 constant functions
belong to $\woo$.)

\begin{proposition}
  \label{propin}
There exists a smooth (complex-valued) 
function $v_0\in \woo(\R)$ such that $v_0(-\infty)=1$, 
$|v_0(x)|\ge\varepsilon_0> 0$ for any $x\in\R$, and the solution $v$ 
of the heat equation with $v(\cdot,0)=v_0$ 
vanishes at infinitely many points $(0,\tau_k)$, with  $\tau_k\to\infty$.
\end{proposition}
\begin{myproof}
  First choose a smooth  real-valued  odd function $w\in \woo(\R)$ 
such that $w(-\infty)=1$ and $w>0$ on $(-\infty,0)$. 
The solution $v_1$ of the heat equation with $v_1(\cdot,0)=w$
has a unique zero at $x=0$ for each $t$.
We shall next find a  smooth real-valued function $z\in \woo(\R)$ 
such that $z(\pm\infty)=0$, 
$z(0)\ne 0$, and the solution   $v_2$ of the heat equation with 
$v_2(\cdot,0)=z$ vanishes at  points $(0,\tau_k)$ with  
$\tau_k\to\infty$. From this the conclusion of the proposition
follows upon setting $v_0=w+iz$.

We first choose sequences $R_k>0$, $\ep_k\in (0,1)$ with the 
following properties:
\begin{itemize}
\item[(a1)] $\sum_{k=1}^\infty\ep_k<\infty$,
\item[(a2)] $R_{k+1}>R_k$ and $R_{k+1}\ep_{k+1}>\sum_{j=1}^k\ep_j(R_j+1)$
 \quad ($k=1,2\dots$).
\end{itemize}
Next, for each $k$ we choose a smooth function $z_k$ such that
\begin{itemize}
\item[(a3)] $z_k\equiv \ep_k$ on $[-R_k,R_k]$, \quad $z_k\equiv 0$ on
  $\R\setminus [-R_k-1,R_k+1]$, 
\item[(a4)]  $0\le z_k\le \ep_k$ and  $|z'_k|\le 2\ep_k$ on $\R$.
\end{itemize}
We will show that if $x_k$ is a suitably chosen sequence, then the function
\begin{equation}
  \label{zdef}
  z(x)=\sum_{k=1}^\infty(-1)^{k+1}z_k(x-x_k)
\end{equation}
has the desired properties. 

The sequence $x_k$ will be constructed so that, in particular,
\begin{equation}
  \label{seq}
  |x_{j}|> |x_{j-1}|+R_j+R_{j-1}+2
\end{equation}
for $j=2,3\dots$. This guarantees that  the  functions
 $z_k(\,\,\cdot-x_k)$ have non-overlapping supports, 
hence, by  (a1), (a3), and  (a4),
 $z$ is a smooth function in $\woo(\R)$ satisfying 
$z(\pm\infty)=0$. In addition to \eqref{seq}, we
need to ensure that the function
\begin{equation}
  \label{heat}
  v_2(0,t)=\frac{1}{\sqrt{4\pi t}}\int_\R
  e^{\frac{-|y|^2}{4t}}z(y)\,dy
\end{equation}
has infinitely many sign changes.

We shall recursively construct sequences $x_k$, $t_k $, $d_k$ 
such that $(-1)^{k+1}x_k$  $t_k$, $d_k$ are nonnegative and
increasing with $k$ 
and the following statement is satisfied for each $k$:
\newcommand{\Ik}{(I$_k$) }
\begin{itemize}
\item[(I$_k$)]   Relations \eqref{seq} hold for
  $j=2,\dots,k$ and, with any choice of
$x_{k+1}, x_{k+2}, \dots$ satisfying \eqref{seq} 
for $j=k+1,k+2,\dots$ and $|x_{k+1}|\ge d_k$, one has
\begin{equation*}
  (-1)^j \int_\R
  e^{\frac{-|y|^2}{4t_j}}z(y)\,dy<0,\quad j=1,\dots, k.
\end{equation*}
\end{itemize}

Take  $x_1=0$, $t_1=1$. It is obvious that (I$_1$) is satisfied
if $d_1$ is sufficiently large. We fix such a $d_1$ satisfying also
$d_1\ge R_1+R_2+2$.

Assume that $x_j$, $t_j$, $d_j$ have been constructed for
$j=1,\dots,k$. To define the next terms assume for definiteness that
$k$ is odd (in case it is  even, the construction is analogous).
Set $x_{k+1}=-d_k$. Assuming $x_{k+2}, x_{k+3},\dots$, are any
numbers satisfying \eqref{seq} for $j=k+2,k+2,\dots$ 
and $|x_{k+2}|\ge d_{k+1}$, with $d_{k+1}$ to be specified below, we
use  (a3), (a4) to estimate
\begin{align}
\int_\R
  e^{\frac{-|y|^2}{4t}}z(y)\,dy
&=\sum_{j=1}^\infty (-1)^{j+1} 
\int_{x_j-R_j-1}^{x_j+R_j+1}
e^{\frac{-|y|^2}{4t}}z_j(y-x_j)\,dy\notag  \\
 &\le 2\sum_{j=1}^k \ep_j(R_j+1)-\ep_{k+1} \int_{x_{k+1}-R_{k+1}}^{x_{k+1}+R_{k+1}}
e^{\frac{-|y|^2}{4t}}\,dy\label{tt1}\\
&~\qquad\qquad\qquad +\int_{\R\setminus[-d_{k+1}+R_{k+2}+1\,,\,d_{k+1}-R_{k+2}-1]}
  e^{\frac{-|y|^2}{4t}}\,dy.
\label{tt2}
\end{align} 
In view of (a2), if  $t=t_{k+1}>t_k$ is large enough, 
the expression in \eqref{tt1} is  negative. Fixing such a $t_{k+1}$
and subsequently choosing a large enough $d_{k+1}>d_k+R_{k+1}+R_k+2$,
we make the whole  expression in \eqref{tt1}, \eqref{tt2} negative. 
Hence (I$_{k+1}$) is satisfied.

It is obvious that with sequences $x_k$ and $t_k$ resulting from the above
construction, the function $z$ has all the desired properties. In
particular, $v_2(0,t)$ has a zero $\tau_k$ in $(t_k,t_{k+1})$ for each $k$. 
\end{myproof}

\begin{Remarks}
  \label{rm-last}{\rm
 (i) It is not difficult to  check that 
with the initial data $z$  constructed above, 
the solution $v_2$ of the heat equation
has a unique zero $x(t)$ for each $t>0$ 
(and $x(t)$ changes sign infinitely many times as $t\to\infty$).

(ii) It is conceivable that if the initial condition $u_0$ is compactly
supported, the solution $u$ of equation~(\ref{B}) given by the Cole-Hopf
transformation will always have only
finitely many singularities. For example, the construction above
cannot be carried out if we demand that $z$ be compactly supported.
This can be seen from the following observation:
{\sl
If $z$ solves the heat equation in $\R\times(0,\infty)$ with compactly 
supported initial data and $t\to z(0,t)$ has infinitely many zeros,
then $z(0,t)=0$ for all $t>0$.}
To see this we recall that Appell's transformation $\tilde z$ of $z$
is defined in $\R\times(-\infty,0)$ by
$z(x,t)=\Gamma(x,t)\tilde z(x/t,-1/t)$. One can see from
formula~\eqref{bounded} applied to $z$ that $\tilde z$ has 
an analytic extension to $\R\times\R$ and
the assumptions on $z$ imply that $\ttt\to\tilde z(0,\ttt)$ has infinitely
many roots accumulating at $0$. Hence $\tilde z(0,\ttt)=0$ for all $\ttt$,
which implies our statement.
}
\end{Remarks}

\vspace{10pt}
\noindent
{\bf Acknowledgments.} The authors would like to thank Luis
Escauriaza for a valuable 
discussion. 


\begin{thebibliography}
\small

\bibitem{Angenent-F:s1}
S.~Angenent and B.~Fiedler, \emph{The dynamics of rotating waves in scalar
  reaction diffusion equations}, Trans. Amer. Math. Soc. \textbf{307} (1988),
  545--568.

\bibitem{birnir}
B.~Birnir, {\emph{An example of blow-up, for the complex KdV equation 
and existence beyond the blow-up.} SIAM J. Appl. Math. \textbf{47} 
(1987), no. 4, 710--725.

\bibitem{Chen:strong}
X.-Y. Chen, \emph{A strong unique continuation theorem for parabolic
  equations}, Math. Ann. \textbf{311} (1998), 603--630.

\bibitem{Friedman:bk-par}
A.~Friedman, \emph{Partial differential equations of parabolic type},
  Prentice-Hall Inc., Englewood Cliffs, N.J., 1964.

\bibitem{Hopf}
E.~Hopf,  \emph{The partial differential equation $u\sb t+uu\sb x=µu\sb {xx}$,} 
Comm.\ Pure Appl.\ Math. \textbf{3}, (1950). 201--230.

\bibitem{Kato}
 T.~Kato, \emph{Strong $L\sp{p}$-solutions of the Navier-Stokes equation in 
 $ \R\sp{m}$, with applications to weak solutions,}  Math. Z. 
 \textbf{187} (1984), no. 4, 471--480. 
 
\bibitem{Koch-Tataru}
 H.~Koch, D.~Tataru, \emph{Well-posedness for the Navier-Stokes equations,} 
 Adv.\ Math. \textbf{157} (2001), no. 1, 22--35.

\bibitem{Sinai}
D.~Li, Y.~Sinai,
\emph {Blow Ups of Complex Solutions of the 3D-Navier-Stokes System,}
arXiv.org preprint,  physics/0610101.

\bibitem{plsv}
P.~Plech\'a\v c, V.~\v Sver\'ak, \emph{ On self-similar singular solutions of the complex 
Ginzburg-Landau equation.}  Comm. Pure Appl. Math. \textbf(54} (2001), 
no. 10, 1215--1242.

\bibitem{senouf1}
D.~Senouf, \emph{Dynamics and condensation of complex singularities for Burgers' equation. I.}
 SIAM J. Math. Anal. \textbf{28} (1997), no. 6, 1457--1489. 

\bibitem{senouf2}
D.~Senouf, R.~Caflisch, N.~Ercolani, 
\emph{Pole dynamics and oscillations for the complex 
Burgers equation in the small-dispersion limit.} Nonlinearity \textbf{9} 
(1996), no. 6, 1671--1702.

\end{thebibliography}
\end{document}